\documentclass[fleqn,10pt]{article}

\usepackage{latexsym,amssymb}

\newtheorem{lemma}{Lemma}

\title{
  A PBW commutator lemma for $U_q[gl(m|n)]$
}

\author{
  David~~De Wit \\
  The University of Queensland, Brisbane, AUSTRALIA \\
  \texttt{ddw@maths.uq.edu.au}
}

\begin{document}

\maketitle

\begin{abstract}
  \noindent
  We present and prove in detail a Poincar\'{e}--Birkhoff--Witt commutator lemma
  for the quantum superalgebra $U_q[gl(m|n)]$.
\end{abstract}


\section{Introduction}

This paper presents and proves in detail a
Poincar\'{e}--Birkhoff--Witt (PBW) commutator lemma for the
quantum superalgebra $U_q[gl(m|n)]$. The lemma itself is not new;
it dates from a 1993 paper of Rui Bin Zhang \cite{Zhang:93} on the
representation theory of $U_q[gl(m|n)]$. However, its previous
incarnation contained several typographical and other minor
errors in its details; and in any case an explicit proof was not
supplied. Here, we correct those errors, and supply detailed
proofs for our claims.

We mention that we use the phrase ``PBW commutator lemma'' to
indicate a result showing commutations sufficient to render any
expression within an algebra into a normal form in a PBW basis;
for more details for our specific case $U_q[gl(m|n)]$, we again
refer the reader to the original work by Zhang.


\section{The structure of $U_q[gl(m|n)]$}
\label{sec:UqglmnFullDescription}

Following Zhang \cite[pp1237-1238]{Zhang:93}, we provide a full
description of $U_q[gl(m|n)]$ in terms of simple generators and
relations. We do so after first introducing the generators and
various divers notations.

Firstly, we define a $\mathbb{Z}_2$ grading $[\cdot]$ on the set
of $gl(m|n)$ indices $\{1, \dots, m+n\}$:
\begin{eqnarray*}
  [a]
  \triangleq
  \left\{
  \begin{array}{lll}
    0 \quad & a \leqslant m\qquad\qquad
            & \mathrm{even~indices}
    \\
    1       & a > m
            & \mathrm{odd~indices},
  \end{array}
  \right.
\end{eqnarray*}
where we use the symbol ``$\triangleq$'' to mean ``is defined as
being''. Throughout, we shall use dummy indices $a,b$, etc.,
where meaningful.

A set of generators for the associative superalgebra
$U_q[gl(m|n)]$ is then:
\begin{eqnarray*}
  \left\{
    K_{a}^{\pm};
    {E^{a}}_{b}
    \; | \;
    1 \leqslant a, b \leqslant m+n,
    \;
    a\neq b
  \right\},
\end{eqnarray*}
where the $K_{a}^{\pm}$ are called ``Cartan generators'' (and of
course we intend ``$\pm1$'' where we write ``$\pm$''), and
${E^{a}}_{b}$ is called a ``raising generator'' if $a<b$ and a
``lowering generator'' if $a>b$. We indeed intend that $K_a$ and
$K_a^{-1}$ are inverses, that is, that we have relations
$K_{a}K_{a}^{-1}=K_{a}^{-1}K_a=\mathrm{Id}$, where $\mathrm{Id}$
is the $U_q[gl(m|n)]$ identity element.

Elements of $U_q[gl(m|n)]$ are then in general weighted sums of
noncommuting products of these generators, where each weight is in
general a rational expression of integer-coefficient Laurent
polynomials in the polynomial variable $q$. Under the phrase
``products of generators'', we include powers of the $K_a$ (see
below).

For various invertible $X$, we will repeatedly use the notation
$\overline{X}\triangleq X^{-1}$; in particular, we set
$\overline{q}\triangleq q^{-1}$. Next, for any index $a$ we shall
write:
\begin{eqnarray*}
  q_a
  \triangleq
  q^{{(-)}^{[a]}},
\end{eqnarray*}
where we have invoked the shorthand ``$(-)$'' for ``$(-1)$''. For
any power $N$, replacing $q$ with $q^N$ immediately shows that
$(q_a)^N=(q^N)_a$, so we may write $q_a^N$ with impunity; in
particular, we will write $\overline{q}_a\equiv q_a^{-1}$.
Further, we will use the following notation:
\begin{eqnarray*}
  \begin{array}{r@{\hspace{5pt}}c@{\hspace{5pt}}lr@{\hspace{5pt}}c@{\hspace{5pt}}l}
    \Delta
    & \triangleq &
    q - \overline{q},
    \qquad
    &
    \Delta_{a}
    & \triangleq &
    q_a-\overline{q}_a
    =
    {(-)}^{[a]}
    (q-\overline{q})
    =
    {(-)}^{[a]}
    \Delta,
    \\
    \overline{\Delta}
    & \triangleq &
    (\Delta)^{-1},
    &
    \overline{\Delta}_{a}
    & \triangleq &
    (\Delta_{a})^{-1}
    =
    {(-)}^{[a]} \overline{\Delta}.
  \end{array}
\end{eqnarray*}

Now, in terms of $q$, an equivalent notation for $K_{a}$ is
$q_a^{{E^a}_a}$. (Here, the exponentiation may be understood in
terms of a power series expansion of the $U[gl(m|n)]$ Cartan
generators ${E^a}_a$. Strictly speaking, we \emph{could} define
these ${E^a}_a$ as the $U_q[gl(m|n)]$ Cartan generators, allowing
them to appear in infinite sums as exponents of $q$, but the $K_a$
notation is more convenient.) Thus, powers $K_{a}^N$ are
meaningful, although we will only deal with
$N\in\frac{1}{2}\mathbb{Z}$ (that is, integer and half-integer
powers). So, we may write $\overline{K}_a\triangleq K_a^{-1}$;
indeed the mapping $q\mapsto \overline{q}$ sends $K_{a}^N$ to
$\overline{K}_a^N$, and as expected, for arbitrary powers $M,N$:
\begin{eqnarray*}
  K_{a}^{M}
  K_{a}^{N}
  =
  K_{a}^{M+N}
  \qquad
  \mathrm{where}
  \qquad
  K_{a}^0
  \equiv
  \mathrm{Id}.
  \qquad
\end{eqnarray*}
Apart from $N\in \mathbb{N}$, powers (i.e. products) of the
non-Cartan generators $({E^{a}}_{b})^N$ for $a\neq b$, are not
meaningful.

The generators inherit a $\mathbb{Z}_2$ grading from the indices:
\begin{eqnarray*}
  [K_{a}]
  \triangleq
  0
  \qquad
  \mathrm{and}
  \qquad
  [{E^a}_b]
  \triangleq
  [a] + [b]
  \quad
  (\mathrm{mod}\;2),
\end{eqnarray*}
so we may also use the terms ``even'' and ``odd'' for generators.
Elements of $U_q[gl(m|n)]$ are said to be \emph{homogeneous} if
they are `linear' combinations of generators of the same grading
or products of other homogeneous elements; the product $XY$ of
homogeneous $X,Y$ has grading $
  [ X Y ]
  \triangleq
  [X] + [Y]
  \;
  (\mathrm{mod}\;2)
$.

Now, the full set of generators includes some redundancy; in that
its elements may be expressed in terms of a subset of them, that
is the following \emph{$U_q[gl(m|n)]$ simple generators}:
\begin{eqnarray*}
  \left\{
    K_{a}^{\pm};
    {E^{a+1}}_{a},
    {E^{a}}_{a+1}
    \; | \;
    1\leqslant a,a+1\leqslant m+n
  \right\};
\end{eqnarray*}
note that there are only two \emph{odd} simple generators:
${E^{m+1}}_{m}$ (lowering) and ${E^{m}}_{m+1}$ (raising). In the
$gl(m|n)$ case, the remaining \emph{nonsimple} (non-Cartan)
generators satisfy the same commutation relations as the simple
generators. However, for $U_q[gl(m|n)]$, the nonsimple generators
are instead recursively defined in terms of weighted sums of
products of simple generators (\cite[p1971,~(3)]{Zhang:92} and
\cite[p1238,~(2)]{Zhang:93}). Writing
$S^a_b\triangleq\mathrm{sign}(a-b)$, the elements of the set of
nonsimple generators $\{{E^{a}}_{b} \; | \; |a-b|>1\}$ may be
defined by:
\begin{eqnarray}
  {E^{a}}_{b}
  \triangleq
  {E^{a}}_{c}
  {E^{c}}_{b}
  -
  q_c^{S^a_b}
  {E^{c}}_{b}
  {E^{a}}_{c},
  \label{eq:UqglmnNonSimpleGeneratorsDefn}
\end{eqnarray}
where we intend $c$ to be an \emph{arbitrary} index strictly
between $a$ and $b$; we do \emph{not} intend a sum here.

Lastly, the \emph{graded commutator} $[\cdot,\cdot]$ is defined
for homogeneous $X,Y$ by:
\begin{eqnarray}
  [ X, Y ]
  \triangleq
  X Y - {(-)}^{[X][Y]} Y X,
  \label{eq:glmnCommutatorBracket}
\end{eqnarray}
and extended by linearity.  As $U_q[gl(m|n)]$ is an
\emph{associative} superalgebra, we have the following useful
identities involving homogeneous elements:
\begin{eqnarray}
  \left.
  \begin{array}{lr@{\hspace{5pt}}c@{\hspace{5pt}}l}
    (\mathrm{a}) &
    [XY,Z]
    & = &
    X[Y,Z]
    +
    {(-)}^{[Y][Z]}
    [X,Z]Y
    \\
    (\mathrm{b}) &
    [X,YZ]
    & = &
    [X,Y]Z
    +
    {(-)}^{[X][Y]}
    Y[X,Z].
  \end{array}
  \hspace{83pt}
  \right\}
  \label{eq:AssociativeSAIdentity}
\end{eqnarray}


\subsection{$U_q[gl(m|n)]$ relations}
\label{sec:Uqglmnrelations}

In terms of the set of simple generators, that is:
\begin{eqnarray*}
  \left\{
    K_{a}^{\pm};
    {E^{a+1}}_{a},
    {E^{a}}_{a+1}
    \; | \;
    1\leqslant a,a+1\leqslant m+n
  \right\},
\end{eqnarray*}
our algebra $U_q[gl(m|n)]$ satisfies the following relations:

\begin{enumerate}
\item
  The Cartan generators commute, that is for $M,N\in\{\pm1\}$:
  \begin{eqnarray}
    K_{a}^{M}
    K_{b}^{N}
    =
    K_{b}^{N}
    K_{a}^{M},
    \qquad
    \qquad
    K_{a} \overline{K}_a
    =
    \mathrm{Id}.
    \label{eq:CartanGeneratorsCommute}
  \end{eqnarray}

\item
  The Cartan generators commute with the simple raising and lowering
  generators in the following manner:
  \begin{eqnarray}
    K_{a}
    {E^b}_{b\pm1}
    =
    q_a^{( \delta^a_b - \delta^a_{b\pm1} )}
    {E^b}_{b\pm1}
    K_{a}.
    \label{eq:SimpleCNCCommutator}
  \end{eqnarray}

\item
  The non-Cartan simple generators satisfy:
  \begin{eqnarray}
    [ {E^a}_{a+1}, {E^{b+1}}_b ]
    =
    \delta^a_b
    \overline{\Delta}_a
    \left(
      K_{a} \overline{K}_{a+1} - \overline{K}_a K_{a+1}
    \right),
    \label{eq:SimpleRLSymmetricCommutator}
  \end{eqnarray}
  and, for $|a-b| > 1$, we have the commutations:
  \begin{eqnarray}
    \hspace{-7pt}
    {E^{a+1}}_{a}
    {E^{b+1}}_{b}
    =
    {E^{b+1}}_{b}
    {E^{a+1}}_{a}
    \quad
    \mathrm{and}
    \quad
    {E^{a}}_{a+1}
    {E^{b}}_{b+1}
    =
    {E^{b}}_{b+1}
    {E^{a}}_{a+1}.
    \label{eq:MoreSimpleRLSymmetricCommutator}
  \end{eqnarray}

\item
  The squares of the odd simple generators are zero:
  \begin{eqnarray}
    {( {E^{m}}_{m+1} )}^2
    =
    {( {E^{m+1}}_{m} )}^2
    =
    0.
    \label{eq:SquaresofOddSimpleGeneratorsareZero}
  \end{eqnarray}

\item
  If neither $m$ nor $n$ is $1$, we have the $U_q[gl(m|n)]$
  \emph{Serre relations} (else if either $m$ or $n$ is $1$, omit them).
  Most succinctly expressed in terms of the
  nonsimple generators, for $a\neq m$, we have:
  \begin{eqnarray}
    \left.
    \begin{array}{cr@{\hspace{5pt}}c@{\hspace{5pt}}l}
      (\mathrm{a}) &
      {E^{a+1}}_{a}
      {E^{a+2}}_{a}
      &=&
      q_{a}
      {E^{a+2}}_{a}
      {E^{a+1}}_{a}
      \\[0.5mm]
      (\mathrm{b}) &
      {E^{a}}_{a+1}
      {E^{a}}_{a+2}
      &=&
      q_{a}
      {E^{a}}_{a+2}
      {E^{a}}_{a+1}
      \\[0.5mm]
      (\mathrm{c}) &
      {E^{a+1}}_{a-1}
      {E^{a+1}}_{a}
      &=&
      q_a
      {E^{a+1}}_{a}
      {E^{a+1}}_{a-1}
      \\[0.5mm]
      (\mathrm{d}) &
      {E^{a-1}}_{a+1}
      {E^{a}}_{a+1}
      &=&
      q_a
      {E^{a}}_{a+1}
      {E^{a-1}}_{a+1},
    \end{array}
    \hspace{70pt}
    \right\}
    \label{eq:SerreRelationsaneqm}
  \end{eqnarray}
  and also:
  \begin{eqnarray*}
    \left[ {E^{m+1}}_{m}, {E^{m+2}}_{m-1} \right]
    =
    \left[ {E^{m}}_{m+1}, {E^{m-1}}_{m+2} \right]
    =
    0.
  \end{eqnarray*}
  The interested reader may use (\ref{eq:UqglmnNonSimpleGeneratorsDefn}) to
  expand these into expressions involving only the simple
  generators; however the results are cumbersome and unedifying.

\end{enumerate}


\subsection{Useful results from the $U_q[gl(m|n)]$ relations}
\label{sec:Uqglmnrelationsuseful}

\begin{enumerate}
\item
  From (\ref{eq:CartanGeneratorsCommute}), it immediately follows that
  all powers of the Cartan generators commute; that is, for any powers
  $M,N\in\frac{1}{2}\mathbb{Z}$:
  \begin{eqnarray}
    K_{a}^{M}
    K_{b}^{N}
    =
    K_{b}^{N}
    K_{a}^{M}.
    \label{eq:PowersofCartanGeneratorsCommute}
  \end{eqnarray}

\item
  Lemma 2 of \cite{DeWit:99c} shows that (\ref{eq:SimpleCNCCommutator})
  may be much strengthened
  to cover all non-Cartan generators and all powers of Cartan generators:
  \begin{eqnarray}
    K_{a}^N
    {E^b}_{c}
    =
    q_a^{N (\delta^a_b-\delta^a_{c})}
    {E^b}_{c}
    K_{a}^N,
    \label{eq:FullCNCCommutator}
  \end{eqnarray}
  that is, where $b, c$ are \emph{any} meaningful indices
  (i.e. even including the case $b=c$),
  and $N\in\frac{1}{2}\mathbb{Z}$ is any power.

\end{enumerate}

The proof of our PBW commutator lemma uses these results, and also
calls on Lemma 1 of \cite{Zhang:93}, which we now cite, with some
slight notational changes and simplifications:
\begin{lemma}
  Where $a<b$, we have the following two results.
  \vspace{\baselineskip}

  \noindent
  Firstly, if $a,b\neq c,c+1$, then:
  \begin{eqnarray}
    \hspace{-30pt}
    \left.
    \begin{array}{cr@{\hspace{5pt}}c@{\hspace{5pt}}l}
      (\mathrm{a}) &
      \left[{E^{a}}_{b}, {E^{c}}_{c+1}\right]
      & = &
      0
      \\
      (\mathrm{b}) &
      \left[{E^{b}}_{a}, {E^{c+1}}_{c}\right]
      & = &
      0.
    \end{array}
    \hspace{206pt}
    \right\}
    \label{eq:Zhang93Lemma1part1}
  \end{eqnarray}
  Secondly, if $a\neq c$ or $b\neq c+1$, then:
  \begin{eqnarray}
    \hspace{-30pt}
    \left.
    \begin{array}{lrcl}
      (\mathrm{a}) &
      \left[{E^{a}}_{b}, {E^{c+1}}_{c}\right]
      &\hspace{-5pt}=\hspace{-5pt}&
      \delta^{c+1}_{b}
      K_{c}
      \overline{K}_{c+1}
      {E^{a}}_{c}
      -
      \delta^{a}_{c}
      {(-)}^{[{E^{c+1}}_{c}]}
      {E^{c+1}}_{b}
      \overline{K}_{c}
      K_{c+1}
      \\
      (\mathrm{b}) &
      \left[{E^{b}}_{a}, {E^{c}}_{c+1}\right]
      &\hspace{-5pt}=\hspace{-5pt}&
      \delta^{c}_{a}
      K_{c}
      \overline{K}_{c+1}
      {E^{b}}_{c+1}
      -
      \delta^{b}_{c+1}
      {(-)}^{[{E^{c}}_{c+1}]}
      {E^{c}}_{a}
      \overline{K}_{c}
      K_{c+1}.
    \end{array}
    \hspace{13pt}
    \right\}
    \label{eq:Zhang93Lemma1part2}
  \end{eqnarray}
\end{lemma}


\subsection{The algebra antiautomorphism $\omega$}

Again following Zhang \cite{Zhang:93}, we introduce an
\emph{ungraded} $U_q[gl(m|n)]$ algebra antiautomorphism $\omega$,
defined for \emph{simple} generators ${E^a}_b$ by:
\begin{eqnarray}
  \omega({E^a}_b)
  \triangleq
  {E^b}_a,
  \qquad
  \omega (K_a)
  \triangleq
  \overline{K}_a,
  \qquad
  \omega (q)
  \triangleq
  \overline{q},
  \label{eq:Defnofomega}
\end{eqnarray}
where by $\omega(q)=\overline{q}$, we intend the more intelligible
$\omega(q \,\mathrm{Id})=\overline{q} \,\mathrm{Id}$. Declaring
$\omega$ to be an ungraded antiautomorphism means that we intend:
\begin{eqnarray}
  \omega(XY)
  =
  \omega(Y)
  \omega(X)
  \qquad
  \mathrm{and}
  \qquad
  \omega(X + Y)
  =
  \omega(X) + \omega(Y);
  \label{eq:omegaXY}
\end{eqnarray}
observe that $\omega$ does indeed preserve grading, that is for
homogeneous $X$, we have $[\omega(X)]=[X]$. Then, for homogeneous
$X,Y$, we have, using (\ref{eq:glmnCommutatorBracket}):
\begin{eqnarray}
  \omega([X,Y])
  =
  [\omega(Y),\omega(X)].
  \label{eq:omegaXcommaY}
\end{eqnarray}

The expression $\omega({E^a}_b)={E^b}_a$ in fact holds for
\emph{all} ${E^a}_b$; the generalisation to nonsimple generators
follows from the application of $\omega$ to their definition in
(\ref{eq:UqglmnNonSimpleGeneratorsDefn}). Moreover, we have
immediately from (\ref{eq:Defnofomega}) the following useful
results:
\begin{eqnarray*}
  \hspace{-5pt}
  \omega (K_a^N)
  =
  \overline{K}_a^N,
  \qquad
  \omega (q^N)
  =
  \overline{q}^N,
  \qquad
  \omega (q_a^N)
  =
  \overline{q}_a^N,
  \qquad
  \omega (\Delta_a)
  =
  - \Delta_a.
\end{eqnarray*}

Zhang goes on to define a set of ``generalised Lusztig
automorphisms'', but we do not require these. In fact, it appears
to be impossible to define them consistently for superalgebras
(as claimed in \cite{Zhang:93}), hence invalidating their use in
the proof of the PBW commutator lemma.

\pagebreak


\section{The PBW commutator lemma}

Using the above machinery, we are now ready to state and prove the
$U_q[gl(m|n)]$ PBW commutator lemma. To whit, we will prove the
following, which is slightly different from the original (Lemma~2
of \cite{Zhang:93}).

\begin{lemma}
  We have the following commutations.
  \vspace{\baselineskip}

  \noindent
  Firstly, \emph{(\ref{eq:SimpleRLSymmetricCommutator})} generalises to the
  case of nonsimple generators, that is:
  \begin{eqnarray}
    \left[{E^{a}}_{b},{E^{b}}_{a}\right]
    \hspace{5pt} = \hspace{5pt}
    \overline{\Delta}_a
    (
      K_{a} \overline{K}_{b}
      -
      \overline{K}_{a} K_{b}
    )
    \hspace{55pt}
    \mathrm{all~}
    a, b.
    \label{eq:TomWaits}
  \end{eqnarray}
  Secondly, where there are three distinct indices, we have:
  \begin{eqnarray}
    \left[{E^{a}}_{c},{E^{c}}_{b}\right]
    & = &
    \left\{
      \begin{array}{ll@{\hspace{60pt}}l}
        (\mathrm{a}) & \overline{K}_{b} K_{c} {E^{a}}_{b} & c<b<a \\
        (\mathrm{b}) & {E^{a}}_{b} K_{a} \overline{K}_{c} & c<a<b \\
        (\mathrm{c}) & {E^{a}}_{b} \overline{K}_{a} K_{c} & b<a<c \\
        (\mathrm{d}) & K_{b} \overline{K}_{c} {E^{a}}_{b} & a<b<c
      \end{array}
      \hspace{30pt}
    \right\}
    \label{eq:BrendanPerry}
    \\
    \left[{E^{c}}_{a},{E^{c}}_{b}\right]
    & = &
    \left[{E^{a}}_{c},{E^{b}}_{c}\right]
    =
    0
    \label{eq:AustralianCrawl}
    \hspace{45pt}
    a<c<b
    \mathrm{~~or~~}
    b<c<a
    \\
    &&
    \hspace{-24mm}
    \left.
    \begin{array}{rcl}
      {E^{c}}_{a} {E^{c}}_{b}
      & = &
      \left\{
        \begin{array}{ll@{\hspace{30pt}}l}
          (\mathrm{a}) & {(-)}^{[{E^{c}}_{b}]} q_{c} {E^{c}}_{b} {E^{c}}_{a}
          & a<b<c \\
          (\mathrm{b}) & {(-)}^{[{E^{c}}_{a}]} q_{c} {E^{c}}_{b} {E^{c}}_{a}
          & c<a<b
        \end{array}
      \right.
      \\
      {E^{a}}_{c} {E^{b}}_{c}
      & = &
      \left\{
        \begin{array}{ll@{\hspace{30pt}}l}
          (\mathrm{c}) & {(-)}^{[{E^{b}}_{c}]} q_{c} {E^{b}}_{c} {E^{a}}_{c}
          & a<b<c \\
          (\mathrm{d}) & {(-)}^{[{E^{a}}_{c}]} q_{c} {E^{b}}_{c} {E^{a}}_{c}
          & c<a<b.
        \end{array}
      \right.
    \end{array}
    \hspace{20pt}
    \right\}
    \label{eq:LeonardCohen}
  \end{eqnarray}
  Thirdly, we describe the situation where there are no common indices,
  where we have $a<b$ and $c<d$.  For $i,j \in \mathbb{N}$,
  let $S(i,j)$ denote the set
  $\{i,i+1,\dots,j\}$. Then, if $S(a,b)$ and $S(c,d)$ are either
  disjoint or one is wholly contained within the other, that is if
  $a<c<d<b$, $a<b<c<d$, $c<a<b<d$ or $c<d<a<b$,
  we have a total of $16$ cases:
  \begin{eqnarray}
    \left[{E^{a}}_{b},{E^{c}}_{d}\right]
    =
    \left[{E^{a}}_{b},{E^{d}}_{c}\right]
    =
    \left[{E^{b}}_{a},{E^{c}}_{d}\right]
    =
    \left[{E^{b}}_{a},{E^{d}}_{c}\right]
    =
    0.
    \label{eq:DrHooks16cases}
  \end{eqnarray}
  More interestingly, if there is some other overlap between the sets
  $S(a,b)$ and $S(c,d)$, that is if $a<c<b<d$ or $c<a<d<b$,
  then we have the $8$ cases:
  \begin{eqnarray}
    & & \hspace{-10mm}
    \left.
    \begin{array}{rcl}
      \left[{E^{a}}_{b},{E^{c}}_{d}\right]
      & = &
      \left\{
        \begin{array}{ll@{\hspace{47pt}}l}
          (\mathrm{a}) & + \Delta_{b} {E^{a}}_{d} {E^{c}}_{b} & a<c<b<d
          \\
          (\mathrm{b}) & - \Delta_{d} {E^{a}}_{d} {E^{c}}_{b} & c<a<d<b
        \end{array}
      \right.
      \\
      \left[{E^{b}}_{a},{E^{d}}_{c}\right]
      & = &
      \left\{
        \begin{array}{ll@{\hspace{47pt}}l}
          (\mathrm{c}) & + \Delta_{b} {E^{d}}_{a} {E^{b}}_{c} & a<c<b<d
          \\
          (\mathrm{d}) & - \Delta_{d} {E^{d}}_{a} {E^{b}}_{c} & c<a<d<b
        \end{array}
     \right.
    \end{array}
    \hspace{8pt}
    \right\}
    \label{eq:RyCooder}
    \\
    & &  \hspace{-10mm}
    \left.
    \begin{array}{rcl}
    \left[{E^{a}}_{b},{E^{d}}_{c}\right]
    & = &
    \left\{
      \begin{array}{ll@{\qquad}l}
        (\mathrm{a}) & - \Delta_{b} \overline{K}_{b} K_{c} {E^{a}}_{c} {E^{d}}_{b} & a<c<b<d \\
        (\mathrm{b}) & + \Delta_{d} {E^{d}}_{b} {E^{a}}_{c} \overline{K}_{a} K_{d} & c<a<d<b
      \end{array}
    \right.
    \\
    \left[{E^{b}}_{a},{E^{c}}_{d}\right]
    & = &
    \left\{
      \begin{array}{ll@{\qquad}l}
        (\mathrm{c}) & - \Delta_{c} {E^{b}}_{d} {E^{c}}_{a} \overline{K}_{c} K_{b} & a<c<b<d \\
        (\mathrm{d}) & + \Delta_{a} \overline{K}_{d} K_{a} {E^{c}}_{a} {E^{b}}_{d} &
        c<a<d<b.
      \end{array}
    \right.
    \end{array}
   \hspace{5pt}
    \right\}
    \label{eq:VishwaMohanBhatt}
  \end{eqnarray}
  \label{lem:PBW}
\end{lemma}

In the above, we disagree with the results published in
\cite{Zhang:93} in several places. Firstly
(\ref{eq:FullCNCCommutator}) shows that
(\ref{eq:BrendanPerry}a,d) are actually equivalent to the
published results:
\begin{eqnarray*}
  \left[{E^{a}}_{c},{E^{c}}_{b}\right]
  =
  \left\{
    \begin{array}{lll}
      (\mathrm{a}) &
      q_{b}
      {E^{a}}_{b}
      K_{c}
      \overline{K}_{b}
      & c<b<a
      \\
      (\mathrm{d}) &
      \overline{q}_{b}
      {E^{a}}_{b}
      K_{b}
      \overline{K}_{c}
      & a<b<c.
    \end{array}
  \right.
\end{eqnarray*}
However, for all the commutators involving no common indices, we
differ in substance.  The published results for
(\ref{eq:RyCooder}) are:
\begin{eqnarray*}
  \left[{E^{a}}_{b},{E^{c}}_{d}\right]
  & \hspace{-5pt} = \hspace{-5pt} &
  +
  \Delta_{b}
  {E^{a}}_{d}
  {E^{c}}_{b}
  \qquad\qquad
  a<c<b<d, \quad c<a<d<b
  \\
  \left[{E^{b}}_{a},{E^{d}}_{c}\right]
  & \hspace{-5pt} = \hspace{-5pt} &
  -
  \Delta_{b}
  {E^{b}}_{c}
  {E^{d}}_{a}
  \qquad\qquad
  a<c<b<d, \quad c<a<d<b,
\end{eqnarray*}
and for (\ref{eq:VishwaMohanBhatt}) are:
\begin{eqnarray*}
  \left[{E^{a}}_{b},{E^{d}}_{c}\right]
  & \hspace{-5pt} = \hspace{-5pt} &
  \left\{
    \begin{array}{ll@{\hspace{20pt}}l}
      (\mathrm{a}) &
      +
      \Delta_{b}
      {E^{d}}_{b}
      {E^{a}}_{c}
      \overline{K}_{b}
      K_{a}
      &
      a<c<b<d
      \\
      (\mathrm{b}) &
      +
      \Delta_{a}
      {E^{a}}_{c}
      {E^{d}}_{b}
      \overline{K}_{a}
      K_{d}
      &
      c<a<d<b.
    \end{array}
  \right.
  \\
  \left[{E^{b}}_{a},{E^{c}}_{d}\right]
  & \hspace{-5pt} = \hspace{-5pt} &
  \left\{
    \begin{array}{ll@{\hspace{20pt}}l}
      (\mathrm{c}) &
      -
      \Delta_{b}
      \overline{K}_{a}
      K_{b}
      {E^{c}}_{a}
      {E^{b}}_{d}
      &
      a<c<b<d
      \\
      (\mathrm{d}) &
      -
      \Delta_{a}
      K_{a}
      \overline{K}_{d}
      {E^{b}}_{d}
      {E^{c}}_{a}
      &
      c<a<d<b.
    \end{array}
  \right.
\end{eqnarray*}
We mention that it was the discovery of errors in
\emph{computations} whilst working on material described in
\cite{DeWit:99c} that led us to check and correct these PBW
results, and consequently rediscover and debug the proof.

\hspace{5mm}

\noindent\textbf{Proof of Lemma \ref{lem:PBW}:}

We prove the components of the lemma in a different order to that
in which we state them. This is to ensure consistency as later
parts of the proof recycle results previously shown.

\begin{itemize}
\item
  \textbf{(\ref{eq:DrHooks16cases})}: These are the $16$ commutators involving
  $a<b$ and $c<d$, with no overlap betwen $S(a,b)$ and $S(c,d)$.

  Firstly, in the cases $a<b<c<d$ and $a<c<d<b$, in evaluating
  $[{E^{a}}_{b},{E^{c}}_{d}]$, we may use
  (\ref{eq:UqglmnNonSimpleGeneratorsDefn}) to recursively expand the
  raising generator ${E^c}_d$ into a sum of products of simple raising
  generators, and then apply (\ref{eq:AssociativeSAIdentity}b) until we
  have a weighted sum of terms all involving commutators of the form
  $[{E^{a}}_{b},{E^{e}}_{e+1}]$, where $a,b\neq e,e+1$, all of which are
  necessarily $0$ by (\ref{eq:Zhang93Lemma1part1}a), thus
  $[{E^{a}}_{b},{E^{c}}_{d}]=0$ for these two cases.

  Secondly, swapping
  $a\leftrightarrow c$ and $b\leftrightarrow d$ in these $2$ cases, and
  rearranging then yields $[{E^{a}}_{b},{E^{c}}_{d}]=0$ for the cases
  $c<d<a<b$ and $c<a<b<d$.

  Thirdly, the $4$ cases $[{E^{a}}_{b},{E^{d}}_{c}]=0$ follow by a similar
  argument, calling on (\ref{eq:Zhang93Lemma1part2}a) rather than
  (\ref{eq:Zhang93Lemma1part1}a).

  Lastly, the remaining $8$ cases $[{E^{b}}_{a},{E^{c}}_{d}]=0$ and
  $[{E^{b}}_{a},{E^{d}}_{c}]=0$ follow by the application of $\omega$ to
  the first $8$ cases, and reversing the commutators.

\pagebreak

\item
  \textbf{(\ref{eq:AustralianCrawl})}: Initially, we show
  (\ref{eq:AustralianCrawl}a), that is for the case $a<c<b$ we show
  $\left[{E^{c}}_{a},{E^{c}}_{b}\right]=0$. If in fact $a=c-1$, then the
  result is already known from (\ref{eq:Zhang93Lemma1part2}a), so we
  assume otherwise, that is we consider the
  case $a<c-1<c<b$:
  \begin{eqnarray*}
    \hspace{-36pt}
    \left[{E^{c}}_{a},{E^{c}}_{b}\right]
    &
      \hspace{-5pt}
      \stackrel{(\ref{eq:UqglmnNonSimpleGeneratorsDefn})}{=}
      \hspace{-5pt}
    &
    \left[
      {E^{c}}_{c-1}
      {E^{c-1}}_{a},
      {E^{c}}_{b}
    \right]
    -
    q_{c-1}
    \left[
      {E^{c-1}}_{a}
      {E^{c}}_{c-1},
      {E^{c}}_{b}
    \right]
    \\
    &
      \hspace{-5pt}
      \stackrel{(\ref{eq:AssociativeSAIdentity}a)}{=}
      \hspace{-5pt}
    &
    {E^{c}}_{c-1}
    \left[
      {E^{c-1}}_{a},
      {E^{c}}_{b}
    \right]
    +
    {(-)}^{[{E^{c-1}}_{a}][{E^{c}}_{b}]}
    \left[
      {E^{c}}_{c-1},
      {E^{c}}_{b}
    \right]
    {E^{c-1}}_{a}
    \\
    & & \hspace{-30pt}
    -
    q_{c-1}
    \left(
      {E^{c-1}}_{a}
      \left[
        {E^{c}}_{c-1},
        {E^{c}}_{b}
      \right]
      +
      {(-)}^{[{E^{c}}_{c-1}][{E^{c}}_{b}]}
      \left[
        {E^{c-1}}_{a},
        {E^{c}}_{b}
      \right]
      {E^{c}}_{c-1}
    \right)
    \\
    &
      \hspace{-5pt}
      \stackrel{(\ref{eq:DrHooks16cases})}{=}
      \hspace{-5pt}
    &
    {(-)}^{[{E^{c-1}}_{a}][{E^{c}}_{b}]}
    \left[
      {E^{c}}_{c-1},
      {E^{c}}_{b}
    \right]
    {E^{c-1}}_{a}
    -
    q_{c-1}
    {E^{c-1}}_{a}
    \left[
      {E^{c}}_{c-1},
      {E^{c}}_{b}
    \right]
    \\
    &
      \hspace{-5pt}
      \stackrel{(\ref{eq:Zhang93Lemma1part2}a)}{=}
      \hspace{-5pt}
    &
    0.
  \end{eqnarray*}

  Swapping $a\leftrightarrow b$ and reversing the commutator then
  yields $\left[{E^{c}}_{a},{E^{c}}_{b}\right]=0$ for the case
  $b<c<a$.  Taking $\omega$ of these two cases yields
  $\left[{E^{a}}_{c},{E^{b}}_{c}\right]=0$ for the cases $a<c<b$ and
  $b<c<a$.

\item
  \textbf{(\ref{eq:TomWaits})}: We show the result
  for $a<b$ using `strong' mathematical induction, that is, we assume
  it true for all $a',b'$ such that $|a'-b'|<|a-b|$, and use this to
  show that it is then necessarily true for our $a,b$. To this end, we
  already know from (\ref{eq:SimpleRLSymmetricCommutator}) that it is true for
  $|a-b|=1$. (If $|a-b|\leqslant 1$, the result is
  already true, indeed trivially so if $a=b$.)
  To whit, where $a<b$, and $b-a>1$, that is $a<b-1<b$, we have:
  \begin{eqnarray}
    \hspace{-35pt}
    \left[ {E^{a}}_{b}, {E^{b}}_{a} \right]
    & \hspace{-5pt} \stackrel{(\ref{eq:UqglmnNonSimpleGeneratorsDefn})}{=} \hspace{-5pt} &
    \left[
      {E^{a}}_{b},
      {E^{b}}_{b-1} {E^{b-1}}_{a}
      -
      q_{b-1}
      {E^{b-1}}_{a} {E^{b}}_{b-1}
    \right]
    \nonumber
    \\
    & \hspace{-5pt} \stackrel{(\ref{eq:AssociativeSAIdentity}b)}{=} \hspace{-5pt} &
    \left[ {E^{a}}_{b}, {E^{b}}_{b-1} \right] {E^{b-1}}_{a}
    \nonumber
    +
    {(-)}^{[{E^a}_b][{E^b}_{b-1}]}
    {E^{b}}_{b-1}
    \left[ {E^{a}}_{b}, {E^{b-1}}_{a} \right]
    \nonumber
    \\
    &&
    \hspace{-65pt}
    -
    q_{b-1}
    \left[ {E^{a}}_{b}, {E^{b-1}}_{a} \right] {E^{b}}_{b-1}
    -
    {(-)}^{[{E^a}_b][{E^a}_{b-1}]}
    q_{b-1}
    {E^{b-1}}_{a}
    \left[ {E^{a}}_{b}, {E^{b}}_{b-1} \right],
    \label{eq:Rishikesh5}
  \end{eqnarray}
  where the factors $[{E^a}_b]\equiv [a]+[b]$ within the parity factors are redundant.
  In (\ref{eq:Rishikesh5}), we thus require the evaluation of the
  commutators $\left[ {E^{a}}_{b}, {E^{b}}_{b-1} \right]$ and
  $\left[ {E^{a}}_{b}, {E^{b-1}}_{a} \right]$.  To this end, we have firstly:
  \begin{eqnarray}
    \left[ {E^{a}}_{b}, {E^{b}}_{b-1} \right]
    \stackrel{(\ref{eq:Zhang93Lemma1part2}a)}{=}
    K_{b-1} \overline{K}_{b} {E^{a}}_{b-1},
    \label{eq:Rishikesh6}
  \end{eqnarray}
  and secondly:
  \begin{eqnarray*}
    \hspace{-26pt}
    \left[ {E^{a}}_{b}, {E^{b-1}}_{a} \right]
    & \stackrel{(\ref{eq:UqglmnNonSimpleGeneratorsDefn})}{=} &
    \left[
      {E^{a}}_{b-1} {E^{b-1}}_{b}
      -
      \overline{q}_{b-1}
      {E^{b-1}}_{b} {E^{a}}_{b-1},
      {E^{b-1}}_{a}
    \right]
    \nonumber
    \\
    & \stackrel{(\ref{eq:AssociativeSAIdentity}a)}{=} &
    {E^{a}}_{b-1} \left[ {E^{b-1}}_{b}, {E^{b-1}}_{a} \right]
    \nonumber
    \\
    &&
    +
    {(-)}^{[{E^{b-1}}_b][{E^{b-1}}_{a}]}
    \left[ {E^{a}}_{b-1}, {E^{b-1}}_{a} \right] {E^{b-1}}_{b}
    \nonumber
    \\
    &&
    -
    \overline{q}_{b-1}
    {E^{b-1}}_{b} \left[ {E^{a}}_{b-1}, {E^{b-1}}_{a} \right]
    \nonumber
    \\
    &&
    -
    \overline{q}_{b-1}
    {(-)}^{[{E^{b-1}}_a]}
    \left[ {E^{b-1}}_{b}, {E^{b-1}}_{a} \right] {E^{a}}_{b-1}
    \\
    & \stackrel{(\ref{eq:AustralianCrawl})}{=} &
    \left[ {E^{a}}_{b-1}, {E^{b-1}}_{a} \right] {E^{b-1}}_{b}
    -
    \overline{q}_{b-1}
    {E^{b-1}}_{b} \left[ {E^{a}}_{b-1}, {E^{b-1}}_{a} \right].
  \end{eqnarray*}
  Using the strong inductive assumption, we then have:
  \begin{eqnarray}
    \hspace{-9mm}
    \left[ {E^{a}}_{b}, {E^{b-1}}_{a} \right]
    & = &
    \overline{\Delta}_{a}
    \left(
      \begin{array}{l}
        \left(
          K_{a} \overline{K}_{b-1} - \overline{K}_{a} K_{b-1}
        \right)
        {E^{b-1}}_{b}
        \\
        \qquad
        -
        \overline{q}_{b-1}
        {E^{b-1}}_{b}
        \left( K_{a} \overline{K}_{b-1} - \overline{K}_{a} K_{b-1}\right)
      \end{array}
    \right)
    \nonumber
    \\
    & \stackrel{(\ref{eq:FullCNCCommutator})}{=} &
    \overline{\Delta}_{a}
    {E^{b-1}}_{b}
    \left(
      \begin{array}{l}
        \overline{q}_{b-1}
        K_{a} \overline{K}_{b-1}
        -
        q_{b-1}
        \overline{K}_{a} K_{b-1}
        \\
        \qquad
        -
        \overline{q}_{b-1}
        K_{a} \overline{K}_{b-1}
        +
        \overline{q}_{b-1}
        \overline{K}_{a} K_{b-1}
      \end{array}
    \right)
    \nonumber
    \\
    & = &
    -
    \overline{\Delta}_{a}
    {E^{b-1}}_{b}
    \overline{K}_{a}
    K_{b-1}
    \left( q_{b-1} - \overline{q}_{b-1} \right)
    \nonumber
    \\
    & = &
    -
    \overline{\Delta}
    {(-)}^{[a]}
    \Delta
    {(-)}^{[b-1]}
    {E^{b-1}}_{b}
    \overline{K}_{a}
    K_{b-1}
    \nonumber
    \\
    & = &
    -
    {(-)}^{[{E^{b-1}}_a]}
    {E^{b-1}}_{b}
    \overline{K}_{a}
    K_{b-1}.
    \label{eq:Rishikesh7}
  \end{eqnarray}
  Now substitute (\ref{eq:Rishikesh6}) and (\ref{eq:Rishikesh7}) into
  (\ref{eq:Rishikesh5}):
  \begin{eqnarray*}
    \hspace{-26pt}
    \left[ {E^{a}}_{b}, {E^{b}}_{a} \right]
    & = &
    K_{b-1} \overline{K}_{b} {E^{a}}_{b-1} {E^{b-1}}_{a}
    \\
    &&
    -
    {(-)}^{[{E^b}_{b-1}]}
    {(-)}^{[{E^{b-1}}_a]}
    {E^{b}}_{b-1} {E^{b-1}}_{b}
    K_{b-1} \overline{K}_{a}
    \\
    &&
    +
    {(-)}^{[{E^{b-1}}_a]}
    q_{b-1}
    {E^{b-1}}_{b} K_{b-1} \overline{K}_{a} {E^{b}}_{b-1}
    \\
    &&
    -
    {(-)}^{[{E^{b-1}}_a]}
    q_{b-1}
    {E^{b-1}}_{a} K_{b-1} \overline{K}_{b} {E^{a}}_{b-1}
    \\
    & = &
    \left(
      {E^{a}}_{b-1} {E^{b-1}}_{a}
      -
      {(-)}^{[{E^{b-1}}_a]}
      {E^{b-1}}_{a} {E^{a}}_{b-1}
    \right)
    K_{b-1} \overline{K}_{b}
    \\
    &&
    \hspace{-30pt}
    -
    {(-)}^{[{E^a}_b]}
    \left(
      {E^{b}}_{b-1} {E^{b-1}}_{b}
      -
      {(-)}^{[{E^{b-1}}_b]}
      {E^{b-1}}_{b} {E^{b}}_{b-1}
    \right)
    K_{b-1} \overline{K}_{a}
    \\
    & \stackrel{(\ref{eq:glmnCommutatorBracket})}{=} &
    \left[ {E^{a}}_{b-1}, {E^{b-1}}_{a} \right]
    K_{b-1} \overline{K}_{b}
    \\
    &&
    -
    {(-)}^{[{E^a}_b]}
    \left[ {E^{b}}_{b-1}, {E^{b-1}}_{b} \right]
    K_{b-1} \overline{K}_{a}
    \\
    & = &
    \overline{\Delta}_{a}
    \left( K_{a} \overline{K}_{b-1} - \overline{K}_{a} K_{b-1} \right)
    K_{b-1} \overline{K}_{b}
    \\
    &&
    -
    {(-)}^{[{E^a}_b]}
    \overline{\Delta}_{b}
    \left( K_{b} \overline{K}_{b-1} - \overline{K}_{b} K_{b-1}  \right)
    K_{b-1} \overline{K}_{a}
    \\
    & = &
    \overline{\Delta}_{a}
    \left(
      K_{a} \overline{K}_{b}
      -
      \overline{K}_{a} K_{b-1}^2 \overline{K}_{b}
      -
      K_{b} \overline{K}_{a}
      +
      \overline{K}_{b} K_{b-1}^2 \overline{K}_{a}
    \right)
    \\
    & = &
    \overline{\Delta}_{a}
    \left( K_{a} \overline{K}_{b} - \overline{K}_{a} K_{b} \right).
  \end{eqnarray*}

  Thus, we have shown (\ref{eq:TomWaits}) for general $a<b$.
  The case $a>b$ then follows by swapping $a\leftrightarrow b$ in
  the above, and rearranging.

\item
  \textbf{(\ref{eq:BrendanPerry})}:
  We first show (\ref{eq:BrendanPerry}a), that is for the case $c<b<a$:
  \begin{eqnarray*}
    \hspace{-20pt}
    \left[{E^{a}}_{c},{E^{c}}_{b}\right]
    &
      \hspace{-5pt}
     \stackrel{(\ref{eq:UqglmnNonSimpleGeneratorsDefn})}{=}
      \hspace{-5pt}
    &
    [{E^{a}}_{b} {E^{b}}_{c}, {E^{c}}_{b}]
    -
    q_{b}
    [{E^{b}}_{c} {E^{a}}_{b}, {E^{c}}_{b}]
    \\
    &
      \hspace{-5pt}
      \stackrel{(\ref{eq:AssociativeSAIdentity}a)}{=}
      \hspace{-5pt}
    &
    {E^{a}}_{b}
    [{E^{b}}_{c},{E^{c}}_{b}]
    +
    {(-)}^{[{E^{b}}_{c}]}
    [{E^{a}}_{b}, {E^{c}}_{b}]
    {E^{b}}_{c}
    \\
    & &
    -
    q_{b}
    {E^{b}}_{c}
    [{E^{a}}_{b}, {E^{c}}_{b}]
    -
    {(-)}^{[{E^{a}}_{b}][{E^{c}}_{b}]}
    q_b
    [{E^{b}}_{c}, {E^{c}}_{b}]
    {E^{a}}_{b}
    \\
    &
      \hspace{-5pt}
      \stackrel{(\ref{eq:AustralianCrawl})}{=}
      \hspace{-5pt}
    &
    {E^{a}}_{b}
    [{E^{b}}_{c},{E^{c}}_{b}]
    -
    q_{b}
    [{E^{b}}_{c}, {E^{c}}_{b}]
    {E^{a}}_{b}
    \\
    &
      \hspace{-5pt}
      \stackrel{(\ref{eq:TomWaits})}{=}
      \hspace{-5pt}
    &
    \overline{\Delta}_{b}
    \left(
      {E^{a}}_{b}
      (K_{b} \overline{K}_{c} - \overline{K}_{b} K_{c})
      -
      q_{b}
      (K_{b} \overline{K}_{c} - \overline{K}_{b} K_{c})
      {E^{a}}_{b}
    \right)
    \\
    &
      \hspace{-5pt}
      \stackrel{(\ref{eq:FullCNCCommutator})}{=}
      \hspace{-5pt}
    &
    \overline{\Delta}_{b}
    \left(
      q_{b}
      K_{b}
      \overline{K}_{c}
      -
      \overline{q}_{b}
      \overline{K}_{b}
      K_{c}
      -
      q_{b}
      K_{b}
      \overline{K}_{c}
      +
      q_{b}
      \overline{K}_{b}
      K_{c}
    \right)
    {E^{a}}_{b}
    \\
    &
      \hspace{-5pt}
      =
      \hspace{-5pt}
    &
    \overline{K}_{b}
    K_{c}
    {E^{a}}_{b}.
  \end{eqnarray*}
  A parallel proof yields (\ref{eq:BrendanPerry}c) for the
  case $b<a<c$:
  \begin{eqnarray*}
    \hspace{-20pt}
    \left[{E^{a}}_{c},{E^{c}}_{b}\right]
    &
      \hspace{-5pt}
      \stackrel{(\ref{eq:UqglmnNonSimpleGeneratorsDefn})}{=}
      \hspace{-5pt}
    &
    [
      {E^{a}}_{c},
      {E^{c}}_{a}
      {E^{a}}_{b}
    ]
    -
    q_{a}
    [
      {E^{a}}_{c},
      {E^{a}}_{b}
      {E^{c}}_{a}
    ]
    \\
    &
      \hspace{-5pt}
      \stackrel{(\ref{eq:AssociativeSAIdentity}b)}{=}
      \hspace{-5pt}
    &
    [{E^{a}}_{c},{E^{c}}_{a}]
    {E^{a}}_{b}
    +
    {(-)}^{[{E^{a}}_{c}]}
    {E^{c}}_{a}
    [{E^{a}}_{c}, {E^{a}}_{b}]
    \\
    & &
    -
    q_{a}
    [{E^{a}}_{c}, {E^{a}}_{b}]
    {E^{c}}_{a}
    -
    {(-)}^{[{E^{a}}_{c}][{E^{a}}_{b}]}
    q_a
    {E^{a}}_{b}
    [{E^{a}}_{c}, {E^{c}}_{a}]
    \\
    &
      \hspace{-5pt}
      \stackrel{(\ref{eq:AustralianCrawl})}{=}
      \hspace{-5pt}
    &
    [{E^{a}}_{c},{E^{c}}_{a}]
    {E^{a}}_{b}
    -
    q_{a}
    {E^{a}}_{b}
    [{E^{a}}_{c}, {E^{c}}_{a}]
    \\
    &
      \hspace{-5pt}
      \stackrel{(\ref{eq:TomWaits})}{=}
      \hspace{-5pt}
    &
    \overline{\Delta}_{a}
    \left(
      (K_{a} \overline{K}_{c} - \overline{K}_{a} K_{c})
      {E^{a}}_{b}
      -
      q_{a}
      {E^{a}}_{b}
      (K_{a} \overline{K}_{c} - \overline{K}_{a} K_{c})
    \right)
    \\
    &
      \hspace{-5pt}
      \stackrel{(\ref{eq:FullCNCCommutator})}{=}
      \hspace{-5pt}
    &
    \overline{\Delta}_{a}
    {E^{a}}_{b}
    \left(
      q_{a}
      K_{a}
      \overline{K}_{c}
      -
      \overline{q}_{a}
      \overline{K}_{a}
      K_{c}
      -
      q_{a}
      K_{a}
      \overline{K}_{c}
      +
      q_{a}
      \overline{K}_{a}
      K_{c}
    \right)
    \\
    &
      \hspace{-5pt}
      =
      \hspace{-5pt}
    &
    {E^{a}}_{b}
    \overline{K}_{a}
    K_{c}.
  \end{eqnarray*}
  Taking $\omega$ of (\ref{eq:BrendanPerry}a) yields:
  \begin{eqnarray*}
    \left[{E^{b}}_{c},{E^{c}}_{a}\right]
    \stackrel{(\ref{eq:omegaXY},\ref{eq:omegaXcommaY})}{=}
    {E^{b}}_{a}
    K_{b}
    \overline{K}_{c}
    \qquad
    c<b<a,
  \end{eqnarray*}
  and swapping $a\leftrightarrow b$ then yields
  (\ref{eq:BrendanPerry}b):
  \begin{eqnarray*}
    \left[{E^{a}}_{c},{E^{c}}_{b}\right]
    =
    {E^{a}}_{b}
    K_{a}
    \overline{K}_{c}
    \qquad
    c<a<b.
  \end{eqnarray*}

  Similarly, taking $\omega$ of (\ref{eq:BrendanPerry}c) yields:
  \begin{eqnarray*}
    \left[{E^{b}}_{c},{E^{c}}_{a}\right]
    \stackrel{(\ref{eq:omegaXY},\ref{eq:omegaXcommaY})}{=}
    K_{a}
    \overline{K}_{c}
    {E^{b}}_{a}
    \qquad
    b<a<c,
  \end{eqnarray*}
  and swapping $a\leftrightarrow b$ then yields
  (\ref{eq:BrendanPerry}d):
  \begin{eqnarray*}
    \left[{E^{a}}_{c},{E^{c}}_{b}\right]
    =
    K_{b}
    \overline{K}_{c}
    {E^{a}}_{b}
    \qquad
    a<b<c.
  \end{eqnarray*}

\item
  \textbf{(\ref{eq:LeonardCohen})}: In a sense, these results are
  really glorified Serre relations. We first prove
  (\ref{eq:LeonardCohen}a), that is for the case $a<b<c$.
  Initially assume that $b\neq c-1$ that is $a<b<c-1<c$. Then we have:
  \begin{eqnarray}
    {E^{c}}_{a}
    {E^{c}}_{b}
    &
      \hspace{-5pt}
      \stackrel{(\ref{eq:UqglmnNonSimpleGeneratorsDefn})}{=}
      \hspace{-5pt}
    &
    {E^{c}}_{a}
    \left(
      {E^{c}}_{c-1}
      {E^{c-1}}_{b}
      -
      q_{c-1}
      {E^{c-1}}_{b}
      {E^{c}}_{c-1}
    \right)
    \nonumber
    \\
    &
      \hspace{-5pt}
      \stackrel{(\ref{eq:DrHooks16cases})}{=}
      \hspace{-5pt}
    &
    {E^{c}}_{a}
    {E^{c}}_{c-1}
    {E^{c-1}}_{b}
    -
    {(-)}^{[{E^{c-1}}_b]}
    q_{c-1}
    {E^{c-1}}_{b}
    {E^{c}}_{a}
    {E^{c}}_{c-1}.
    \label{eq:Rishikesh1}
  \end{eqnarray}
  Thus, we must investigate ${E^{c}}_{a} {E^{c}}_{c-1}$. To this end,
  observe that our assumption that $b\neq c-1$ means that we have already
  assumed that $a \neq c-2$, that is, that we safely have $a<c-2<c-1<c$,
  hence:
    \begin{eqnarray}
    {E^{c}}_{a}
    {E^{c}}_{c-1}
    & \stackrel{(\ref{eq:UqglmnNonSimpleGeneratorsDefn})}{=} &
    \left(
      {E^{c}}_{c-2}
      {E^{c-2}}_{a}
      -
      q_{c-2}
      {E^{c-2}}_{a}
      {E^{c}}_{c-2}
    \right)
    {E^{c}}_{c-1}
    \nonumber
    \\
    & \stackrel{(\ref{eq:DrHooks16cases})}{=} &
    {E^{c}}_{c-2}
    {E^{c}}_{c-1}
    {E^{c-2}}_{a}
    -
    q_{c-2}
    {E^{c-2}}_{a}
    {E^{c}}_{c-2}
    {E^{c}}_{c-1}.
    \label{eq:Rishikesh2}
  \end{eqnarray}
  So now, we must investigate ${E^{c}}_{c-2} {E^{c}}_{c-1}$, and this
  falls into two cases. In the general case, if $c \neq m+1$,
  the Serre relation of (\ref{eq:SerreRelationsaneqm}c) gives us:
  $
    {E^{c}}_{c-2}
    {E^{c}}_{c-1}
    =
    q_{c-1}
    {E^{c}}_{c-1}
    {E^{c}}_{c-2}
  $.
  On the other hand, if $c=m+1$, then we have:
  \begin{eqnarray*}
    \hspace{-20pt}
    {E^{m+1}}_{m-1} {E^{m+1}}_{m}
    & \stackrel{(\ref{eq:UqglmnNonSimpleGeneratorsDefn})}{=} &
    \left(
      {E^{m+1}}_{m} {E^{m}}_{m-1}
      -
      q_{m}
      {E^{m}}_{m-1} {E^{m+1}}_{m}
    \right)
    {E^{m+1}}_{m}
    \\
    \hspace{-20pt}
    & \stackrel{(\ref{eq:SquaresofOddSimpleGeneratorsareZero})}{=} &
    {E^{m+1}}_{m} {E^{m}}_{m-1} {E^{m+1}}_{m}
    \\
    \hspace{-20pt}
    {E^{m+1}}_{m} {E^{m+1}}_{m-1}
    & \stackrel{(\ref{eq:UqglmnNonSimpleGeneratorsDefn})}{=} &
    {E^{m+1}}_{m}
    \left(
      {E^{m+1}}_{m} {E^{m}}_{m-1}
      -
      q_{m}
      {E^{m}}_{m-1} {E^{m+1}}_{m}
    \right)
    \\
    \hspace{-20pt}
    & \stackrel{(\ref{eq:SquaresofOddSimpleGeneratorsareZero})}{=} &
    -
    q_{m}
    {E^{m+1}}_{m} {E^{m}}_{m-1} {E^{m+1}}_{m},
  \end{eqnarray*}
  hence
  $
    {E^{m+1}}_{m-1}
    {E^{m+1}}_{m}
    =
    -
    \overline{q}_{m}
    {E^{m+1}}_{m}
    {E^{m+1}}_{m-1}
  $. Taken together, we have for \emph{any} $c$:
  \begin{equation}
    {E^{c}}_{c-2}
    {E^{c}}_{c-1}
    =
    {(-)}^{[{E^c}_{c-1}]}
    q_{c}
    {E^{c}}_{c-1}
    {E^{c}}_{c-2}.
    \label{eq:Rishikesh3}
  \end{equation}
  Installing (\ref{eq:Rishikesh3}) into (\ref{eq:Rishikesh2}), we have:
  \begin{eqnarray}
    \hspace{-35pt}
    {E^{c}}_{a}
    {E^{c}}_{c-1}
    & \hspace{-5pt} = \hspace{-5pt} &
    {(-)}^{[{E^c}_{c-1}]}
    q_{c}
    \left(
      {E^{c}}_{c-1}
      {E^{c}}_{c-2}
      {E^{c-2}}_{a}
      -
      q_{c-2}
      {E^{c-2}}_{a}
      {E^{c}}_{c-1}
      {E^{c}}_{c-2}
     \right)
    \nonumber
    \\
    & \hspace{-5pt} \stackrel{(\ref{eq:DrHooks16cases})}{=} \hspace{-5pt} &
    {(-)}^{[{E^c}_{c-1}]}
    q_{c}
    {E^{c}}_{c-1}
    \left(
      {E^{c}}_{c-2}
      {E^{c-2}}_{a}
      -
      q_{c-2}
      {E^{c-2}}_{a}
      {E^{c}}_{c-2}
    \right)
    \nonumber
    \\
    & \hspace{-5pt} \stackrel{(\ref{eq:UqglmnNonSimpleGeneratorsDefn})}{=} \hspace{-5pt} &
    {(-)}^{[{E^c}_{c-1}]}
    q_{c}
    {E^{c}}_{c-1}
    {E^{c}}_{a}.
    \label{eq:Rishikesh4}
  \end{eqnarray}
  Installing (\ref{eq:Rishikesh4}) into (\ref{eq:Rishikesh1}), we obtain
  the required (\ref{eq:LeonardCohen}a) for the special case
  $a<b<c-1<c$:
  \begin{eqnarray*}
    \hspace{-36pt}
    {E^{c}}_{a}
    {E^{c}}_{b}
    & \hspace{-10pt} = \hspace{-10pt} &
    {(-)}^{[{E^c}_{c-1}]}
    q_{c}
    \left(
      {E^{c}}_{c-1}
      {E^{c}}_{a}
      {E^{c-1}}_{b}
      -
      {(-)}^{[{E^{c-1}}_b]}
      q_{c-1}
      {E^{c-1}}_{b}
      {E^{c}}_{c-1}
      {E^{c}}_{a}
    \right)
    \\
    & \hspace{-10pt} \stackrel{(\ref{eq:DrHooks16cases})}{=} \hspace{-10pt} &
    {(-)}^{[{E^c}_{c-1}]}
    {(-)}^{[{E^{c-1}}_b]}
    q_{c}
    \left(
      {E^{c}}_{c-1}
      {E^{c-1}}_{b}
      -
      q_{c-1}
      {E^{c-1}}_{b}
      {E^{c}}_{c-1}
    \right)
    {E^{c}}_{a}
    \\
    & \hspace{-10pt} \stackrel{(\ref{eq:UqglmnNonSimpleGeneratorsDefn})}{=} \hspace{-10pt} &
    {(-)}^{[{E^c}_b]}
    q_{c}
    {E^{c}}_{b}
    {E^{c}}_{a}.
  \end{eqnarray*}
  If in fact $b=c-1$, then if also $a \neq c-2$, then
  (\ref{eq:Rishikesh4}) covers our result, and if $a=c-2$, then
  (\ref{eq:Rishikesh3}) covers it. Together, we have (\ref{eq:LeonardCohen}a) for all $a<b<c$.

  A parallel proof covers (\ref{eq:LeonardCohen}b), that is, the case
  $c<a<b$; but we omit this. Before proceeding, we condense our notation.
  We have:
  \begin{eqnarray*}
    {E^{c}}_{a}
    {E^{c}}_{b}
    =
    \left\{
      \begin{array}{l@{\qquad}l}
        {(-)}^{[{E^{c}}_{b}]} q_{c} {E^{c}}_{b} {E^{c}}_{a} & a<b<c \\
        {(-)}^{[{E^{c}}_{a}]} q_{c} {E^{c}}_{b} {E^{c}}_{a} & c<a<b.
      \end{array}
    \right.
  \end{eqnarray*}

\pagebreak

  Combining these two results, we may write, for $a<b$:
  \begin{eqnarray}
    {E^{c}}_{a}
    {E^{c}}_{b}
    =
    {(-)}^{[{E^{c}}_{z(a,b,c)}]}
    q_c
    {E^{c}}_{b}
    {E^{c}}_{a}
    \qquad
    \mathrm{if~}
    z(a,b,c) \neq c,
    \label{eq:JJCale}
  \end{eqnarray}
  where $z(a,b,c)$ is a little function which picks out the median element of the set of natural
  numbers $\{a,b,c\}$.
  Applying $\omega$ to (\ref{eq:JJCale}) and cross multiplying yields:
  \begin{eqnarray*}
    {E^{a}}_{c}
    {E^{b}}_{c}
    \stackrel{(\ref{eq:omegaXY})}{=}
    {(-)}^{[{E^{z(a,b,c)}}_c]}
    q_c
    {E^{b}}_{c}
    {E^{a}}_{c}
    \qquad
    \mathrm{if~}
    z(a,b,c) \neq c,
  \end{eqnarray*}
  which is immediately seen to cover (\ref{eq:LeonardCohen}c,d):
  \begin{eqnarray*}
    {E^{a}}_{c}
    {E^{b}}_{c}
    =
    \left\{
      \begin{array}{l@{\qquad}l}
        {(-)}^{[{E^{b}}_{c}]} q_{c} {E^{b}}_{c} {E^{a}}_{c} & a<b<c
        \\
        {(-)}^{[{E^{a}}_{c}]} q_{c} {E^{b}}_{c} {E^{a}}_{c} & c<a<b.
      \end{array}
    \right.
  \end{eqnarray*}

\item
  \textbf{(\ref{eq:RyCooder})}: Beginning with the case
  $a<c<b<d$, we have:
  \begin{eqnarray*}
    \hspace{-40pt}
    \left[ {E^{a}}_{b}, {E^{c}}_{d} \right]
    & \hspace{-10pt} \stackrel{(\ref{eq:glmnCommutatorBracket})}{=} \hspace{-10pt} &
    {E^{a}}_{b} {E^{c}}_{d}
    -
    {(-)}^{[{E^a}_b][{E^c}_d]}
    {E^{c}}_{d} {E^{a}}_{b}
    \\
    & \hspace{-10pt} \stackrel{(\ref{eq:UqglmnNonSimpleGeneratorsDefn})}{=} \hspace{-10pt} &
    {E^{a}}_{b}
    \left(
      {E^{c}}_{b} {E^{b}}_{d}
      -
      \overline{q}_{b} {E^{b}}_{d} {E^{c}}_{b}
    \right)
    \hspace{-1pt}
    -
    \hspace{-1pt}
    {(-)}^{[{E^c}_b]}
    \hspace{-3pt}
    \left(
      {E^{c}}_{b} {E^{b}}_{d}
      -
      \overline{q}_{b} {E^{b}}_{d} {E^{c}}_{b}
    \right)
    {E^{a}}_{b}
    \\
    & \hspace{-10pt} = \hspace{-10pt} &
    \left(
      {E^{a}}_{b} {E^{c}}_{b} {E^{b}}_{d}
      -
      {(-)}^{[{E^c}_b]}
      {E^{c}}_{b} {E^{b}}_{d} {E^{a}}_{b}
    \right)
    \nonumber
    \\
    &&
    \qquad
    -
    \overline{q}_{b}
    \left(
      {E^{a}}_{b} {E^{b}}_{d} {E^{c}}_{b}
      -
      {(-)}^{[{E^c}_b]}
      {E^{b}}_{d} {E^{c}}_{b} {E^{a}}_{b}
    \right).
  \end{eqnarray*}
  Now, for $a<c<b$, by (\ref{eq:LeonardCohen}c), we have
  $
    {E^{a}}_{b} {E^{c}}_{b}
    =
    {(-)}^{[{E^c}_b]}
    q_{b} {E^{c}}_{b} {E^{a}}_{b}
  $.
  Installing this, we quickly obtain (\ref{eq:RyCooder}a):
  \begin{eqnarray*}
    \hspace{-40pt}
    \left[ {E^{a}}_{b}, {E^{c}}_{d} \right]
    & \hspace{-10pt} = \hspace{-10pt} &
    {(-)}^{[{E^c}_b]}
    \hspace{-1pt}
    {E^{c}}_{b}
    \hspace{-1pt}
    \left(
      \hspace{-1pt}
      q_{b}
      {E^{a}}_{b} {E^{b}}_{d}
      -
      {E^{b}}_{d} {E^{a}}_{b}
      \hspace{-1pt}
    \right)
    \hspace{-1pt}
    -
    \hspace{-1pt}
    \overline{q}_{b}
    \hspace{-1pt}
    \left(
      \hspace{-1pt}
      {E^{a}}_{b} {E^{b}}_{d}
      -
      \overline{q}_{b}
      {E^{b}}_{d} {E^{a}}_{b}
      \hspace{-1pt}
    \right)
    \hspace{-1pt}
    {E^{c}}_{b}
    \\
    & \hspace{-10pt} \stackrel{(\ref{eq:UqglmnNonSimpleGeneratorsDefn})}{=} \hspace{-10pt} &
    {(-)}^{[{E^c}_b]}
     q_{b}
    {E^{c}}_{b}  {E^{a}}_{d}
    -
    \overline{q}_{b} {E^{a}}_{d} {E^{c}}_{b}
    \\
    & \hspace{-10pt} \stackrel{(\ref{eq:DrHooks16cases})}{=} \hspace{-10pt} &
    {E^{a}}_{d} {E^{c}}_{b} \left( q_{b} - \overline{q}_{b} \right)
    \\
    & \hspace{-10pt} = \hspace{-10pt} &
    \Delta_{b} {E^{a}}_{d} {E^{c}}_{b}.
  \end{eqnarray*}

  Swapping $a\leftrightarrow c$ and
  $b\leftrightarrow d$ in (\ref{eq:RyCooder}a) then yields:
  \begin{eqnarray}
    \left[{E^{c}}_{d}, {E^{a}}_{b}\right]
    =
    \Delta_{d} {E^{c}}_{b} {E^{a}}_{d}
    \qquad
    c<a<d<b.
    \label{eq:Cake}
  \end{eqnarray}
  Reversing both the commutator and the RHS product yields:
  \begin{eqnarray*}
    -
    {(-)}^{[{E^{c}}_{d}][{E^{a}}_{b}]}
    \left[{E^{a}}_{b}, {E^{c}}_{d}\right]
    \stackrel{(\ref{eq:DrHooks16cases})}{=}
    {(-)}^{[{E^{c}}_{b}][{E^{a}}_{d}]}
    \Delta_{d}
    {E^{a}}_{d}
    {E^{c}}_{b},
  \end{eqnarray*}
  but for $c<a<d<b$, in fact
  $[{E^{c}}_{d}][{E^{a}}_{b}]=[{E^{c}}_{b}][{E^{a}}_{d}]=[{E^a}_d]$, yielding
  (\ref{eq:RyCooder}b):
  \begin{eqnarray*}
    \left[{E^{a}}_{b}, {E^{c}}_{d}\right]
    =
    -
    \Delta_{d}
    {E^{a}}_{d}
    {E^{c}}_{b}
    \qquad
    c<a<d<b.
  \end{eqnarray*}

  Next, applying $\omega$ to (\ref{eq:RyCooder}a) yields:
  \begin{eqnarray*}
    \left[{E^{d}}_{c},{E^{b}}_{a}\right]
    \stackrel{(\ref{eq:omegaXY},\ref{eq:omegaXcommaY})}{=}
    -
    \Delta_{b}
    {E^{b}}_{c}
    {E^{d}}_{a}
    \qquad
    a<c<b<d.
  \end{eqnarray*}
  Reversing both the commutator and the RHS product yields
  (\ref{eq:RyCooder}c):
  \begin{eqnarray*}
    \left[{E^{b}}_{a},{E^{d}}_{c}\right]
    \stackrel{(\ref{eq:DrHooks16cases})}{=}
    \Delta_{b}
    {E^{d}}_{a}
    {E^{b}}_{c}
    \qquad
    a<c<b<d.
  \end{eqnarray*}

  Lastly, applying $\omega$ to
  (\ref{eq:Cake}) yields (\ref{eq:RyCooder}d):
  \begin{eqnarray*}
    \left[{E^{b}}_{a}, {E^{d}}_{c}\right]
    \stackrel{(\ref{eq:omegaXY},\ref{eq:omegaXcommaY})}{=}
    -
    \Delta_{d}
    {E^{d}}_{a}
    {E^{b}}_{c}
    \qquad
    c<a<d<b.
  \end{eqnarray*}

\item
  \textbf{(\ref{eq:VishwaMohanBhatt})}: We first show
  (\ref{eq:VishwaMohanBhatt}a), that is for the case $a<c<b<d$.
  We have:
  \begin{eqnarray*}
    \left[{E^{a}}_{b},{E^{d}}_{c}\right]
    &
      \hspace{-5pt}
     \stackrel{(\ref{eq:UqglmnNonSimpleGeneratorsDefn})}{=}
      \hspace{-5pt}
    &
    \left[
      {E^{a}}_{b},
      {E^{d}}_{b}
      {E^{b}}_{c}
    \right]
    -
    q_{b}
    \left[
      {E^{a}}_{b},
      {E^{b}}_{c}
      {E^{d}}_{b}
    \right]
    \\
    &
      \hspace{-5pt}
      \stackrel{(\ref{eq:AssociativeSAIdentity}b)}{=}
      \hspace{-5pt}
    &
    \left[
      {E^{a}}_{b},
      {E^{d}}_{b}
    \right]
    {E^{b}}_{c}
    +
    {(-)}^{[{E^{d}}_{b}][{E^{a}}_{b}]}
    {E^{d}}_{b}
    \left[
      {E^{a}}_{b},
      {E^{b}}_{c}
    \right]
    \\
    & &
    \hspace{-15pt}
    -
    q_{b}
    \left(
      \left[
        {E^{a}}_{b},
        {E^{b}}_{c}
      \right]
      {E^{d}}_{b}
      +
      {(-)}^{[{E^{a}}_{b}][{E^{b}}_{c}]}
      {E^{b}}_{c}
      \left[
        {E^{a}}_{b},
        {E^{d}}_{b}
      \right]
    \right)
    \\
    &
      \hspace{-5pt}
      \stackrel{(\ref{eq:AustralianCrawl})}{=}
      \hspace{-5pt}
    &
    {E^{d}}_{b}
    \left[{E^{a}}_{b}, {E^{b}}_{c}\right]
    -
    q_{b}
    \left[{E^{a}}_{b}, {E^{b}}_{c}\right]
    {E^{d}}_{b}
    \\
    &
      \hspace{-5pt}
      \stackrel{(\ref{eq:BrendanPerry}d)}{=}
      \hspace{-5pt}
    &
    {E^{d}}_{b}
    K_{c}
    \overline{K}_{b}
    {E^{a}}_{c}
    -
    q_{b}
    K_{c}
    \overline{K}_{b}
    {E^{a}}_{c}
    {E^{d}}_{b}
    \\
    &
      \hspace{-5pt}
    \stackrel{(\ref{eq:FullCNCCommutator},\ref{eq:DrHooks16cases})}{=}
      \hspace{-5pt}
    &
    -
    \Delta_{b}
    \overline{K}_{b}
    K_{c}
    {E^{a}}_{c}
    {E^{d}}_{b}.
  \end{eqnarray*}

  Applying $\omega$ to (\ref{eq:VishwaMohanBhatt}a) yields:
  \begin{eqnarray}
    \left[{E^{c}}_{d},{E^{b}}_{a}\right]
    \stackrel{(\ref{eq:omegaXY},\ref{eq:omegaXcommaY})}{=}
    \Delta_{b}
    {E^{b}}_{d}
    {E^{c}}_{a}
    \overline{K}_{c}
    K_{b}
    \qquad
    a<c<b<d,
    \label{eq:omegaVishwaMohanBhatta}
  \end{eqnarray}
  and swapping $a\leftrightarrow c$ and $b\leftrightarrow d$ then
  yields (\ref{eq:VishwaMohanBhatt}b):
  \begin{eqnarray*}
    \left[{E^{a}}_{b},{E^{d}}_{c}\right]
    =
    \Delta_{d}
    {E^{d}}_{b}
    {E^{a}}_{c}
    \overline{K}_{a}
    K_{d}
    \qquad
    c<a<d<b.
  \end{eqnarray*}

 Next, reversing the commutator in (\ref{eq:omegaVishwaMohanBhatta}) yields:
  \begin{eqnarray*}
    \left[{E^{b}}_{a},{E^{c}}_{d}\right]
    \stackrel{(\ref{eq:omegaXcommaY})}{=}
    -
    {(-)}^{[{E^{b}}_{a}][{E^{c}}_{d}]}
    \Delta_{b}
    {E^{b}}_{d}
    {E^{c}}_{a}
    \overline{K}_{c}
    K_{b}
  \end{eqnarray*}
  However, for the case $a<c<b<d$, we have
  $[{E^{b}}_{a}][{E^{c}}_{d}]=[{E^{b}}_{c}]$, thus:
  $
    {(-)}^{[{E^{b}}_{a}][{E^{c}}_{d}]}
    \Delta_{b}
    =
    {(-)}^{[{E^{b}}_{c}]}
    {(-)}^{[b]}
    \Delta
    =
    {(-)}^{[c]}
    \Delta
    =
    \Delta_{c}
  $,
  yielding (\ref{eq:VishwaMohanBhatt}c):
  \begin{eqnarray*}
    \left[{E^{b}}_{a},{E^{c}}_{d}\right]
    =
    -
    \Delta_{c}
    {E^{b}}_{d}
    {E^{c}}_{a}
    \overline{K}_{c}
    K_{b}
    \qquad
    a<c<b<d.
  \end{eqnarray*}

  Lastly, applying $\omega$ to (\ref{eq:VishwaMohanBhatt}c) yields:
  \begin{eqnarray*}
    \left[{E^{d}}_{c},{E^{a}}_{b}\right]
    \stackrel{(\ref{eq:omegaXY},\ref{eq:omegaXcommaY})}{=}
    \Delta_{c}
    \overline{K}_{b}
    K_{c}
    {E^{a}}_{c}
    {E^{d}}_{b}
    \qquad
    a<c<b<d,
  \end{eqnarray*}
  and then swapping $a\leftrightarrow c$ and $b\leftrightarrow d$
  yields (\ref{eq:VishwaMohanBhatt}d):
  \begin{eqnarray*}
    \left[{E^{b}}_{a},{E^{c}}_{d}\right]
    =
    \Delta_{a}
    \overline{K}_{d}
    K_{a}
    {E^{c}}_{a}
    {E^{b}}_{d}
    \qquad
    c<a<d<b.
  \end{eqnarray*}
\hfill$\Box$
\end{itemize}


\section{Discussion}

Of some interest is that we may use our PBW commutator lemma to
show that (\ref{eq:SquaresofOddSimpleGeneratorsareZero}) in fact
generalises to the nonsimple odd generators, that is:
\begin{eqnarray*}
  {( {E^{a}}_{b} )}^2
  =
  0,
\end{eqnarray*}
for any indices $a,b$ such that $[a] \neq [b]$. The proof of this
statement is left as an (easy) exercise involving
(\ref{eq:LeonardCohen}).

Now that it is established, we may concentrate the notation of our
lemma -- this is useful for encoding purposes.

\begin{itemize}
\item
  The entirety of (\ref{eq:AustralianCrawl}) and
  (\ref{eq:LeonardCohen}) may be summarised by:
  \begin{eqnarray*}
    \hspace{-30pt}
    {E^{a}}_{c}
    {E^{b}}_{c}
    =
    \kappa
    {E^{b}}_{c}
    {E^{a}}_{c}
    \quad
    \mathrm{and}
    \quad
    {E^{c}}_{a}
    {E^{c}}_{b}
    =
    \kappa
    {E^{c}}_{b}
    {E^{c}}_{a},
    \qquad
    \mathrm{any~}
    a \neq b \neq c,
  \end{eqnarray*}
  where:
  \begin{eqnarray*}
    \kappa
    \triangleq
    \left\{
    \begin{array}{l@{\qquad}l}
      1
      &
      \mathrm{if~} z(a,b,c) = c
      \\
      {(-)}^{[{E^{z(a,b,c)}}_c]}
      \overline{q}_{c}^{S^a_b} & \mathrm{else},
    \end{array}
    \right.
  \end{eqnarray*}
  and where $z(a,b,c)$ is our little function which picks out the
  median element of the set of $3$ distinct natural numbers $\{a,b,c\}$. (The $1$ factor follows as
  $[{E^{a}}_{c}][{E^{b}}_{c}]=0$ for $c$ strictly between $a$ and
  $b$.)

\item
  The entirety of (\ref{eq:DrHooks16cases}) to
  (\ref{eq:VishwaMohanBhatt}) may be summarised by:
  \begin{eqnarray*}
    \left[{E^{a}}_{b},{E^{c}}_{d}\right]
    =
    \left\{
      \begin{array}{l@{\hspace{0pt}}l@{\hspace{20pt}}l}
        +
        \Delta_{b} &
        {E^{a}}_{d}
        {E^{c}}_{b}
        &
        a<c<b<d
        \\
        -
        \Delta_{d} &
        {E^{a}}_{d}
        {E^{c}}_{b}
        &
        c<a<d<b
        \\
        +
        \Delta_{a} &
        {E^{c}}_{b}
        {E^{a}}_{d}
        &
        b<d<a<c
        \\
        -
        \Delta_{c} &
        {E^{c}}_{b}
        {E^{a}}_{d}
        &
        d<b<c<a
        \\
        -
        \Delta_{b} &
        \overline{K}_{b}
        K_{d}
        {E^{a}}_{d}
        {E^{c}}_{b}
        &
        a<d<b<c
        \\
        +
        \Delta_{c} &
        {E^{c}}_{b}
        {E^{a}}_{d}
        \overline{K}_{a}
        K_{c}
        &
        d<a<c<b
        \\
        -
        \Delta_{c} &
        {E^{a}}_{d}
        {E^{c}}_{b}
        \overline{K}_{c}
        K_{a}
        &
        b<c<a<d
        \\
        +
        \Delta_{b} &
        \overline{K}_{d}
        K_{b}
        {E^{c}}_{b}
        {E^{a}}_{d}
        &
        c<b<d<a
        \\
        & 0 &
        a\neq b\neq c\neq d \mathrm{~~else}.
      \end{array}
    \right.
  \end{eqnarray*}
\end{itemize}

Finally, we mention that the consistency (if not the veracity) of
our lemma is also supported by extensive computer tests using
\textsc{Mathematica}. By this, we mean that we confirm that:
\begin{eqnarray}
  \mathtt{NormalOrder}(XY)
  =
  \mathtt{NormalOrder}(\mathtt{ExpandNS}(XY)),
  \label{eq:NOTest}
\end{eqnarray}
for a range of $U_q[gl(m|n)]$ nonsimple generators $X,Y$, where
$\mathtt{NormalOrder}(X)$ is a function which renders $X$ in a
normal form, and $\mathtt{ExpandNS}(X)$ is a function which
recursively expands all nonsimple generators in $X$, using
(\ref{eq:UqglmnNonSimpleGeneratorsDefn}).

To be more specific, let the `\emph{height}' of generator
$X\equiv {E^{a}}_{b}$ be $|a-b|$; this is a measure of its
`distance' from simplicity. For $U_q[gl(m|n)]$, it varies from
$0$ (for Cartan generators), to $1$ (for simple non-Cartan
generators); and then for the nonsimple generators from a minimum
of $2$ to a maximum of $m+n-1$ for the `maximally nonsimple'
${E^{m+n}}_{1}$ and ${E^{1}}_{m+n}$.

Then, we confirm that our code satisfies (\ref{eq:NOTest}), for
all $U_q[gl(m|n)]$ generators $X,Y$ of height at most $m+n-1$ for
all $m,n$ such that $m+n\leqslant 5$; at most $3$ for
$m+n\leqslant 10$; and at most $2$ for $m+n\leqslant 18$ (sheer
bloody-mindedness!). The computational expense in performing
these checks rises at least exponentially with height, so we have
to abandon our calculations at this point.  However, our results
do amount to a `complete' consistency check of our lemma, for all
$U_q[gl(m|n)]$ such that $m+n\leqslant 5$.


\section*{Acknowledgements}

I am grateful to Rui Bin Zhang and David McAnally for advice on
issues related to this proof.  The bulk of the work for this
paper was carried out in the year 2000, whilst I was a
Postdoctoral Fellow at Kyoto University in Japan (funded by a
Postdoctoral Fellowship for Foreign Researchers \# P99703,
provided by the Japan Society for the Promotion of Science).
D\={o}mo arigat\={o} gozaimashita! The last algebraic details were
worked out in October 2001, whilst my mind was pacified somewhat
by doing Hatha Yoga at Sri Ved Niketan Ashram at Rishikesh, India.
Dhanyabad! Editing was later performed in March 2002, whilst I was
a guest of Ebadollah Mahmoodian of Sharif University of Technology
in Tehran, Iran. Motashakkeram! Last-minute edits were made in
July 2002, whilst I was a guest of Karin van Dijk in Deventer, The
Netherlands. Dank u wel!

\pagebreak


\bibliographystyle{plain}
\bibliography{DeWit2000}

\end{document}